\documentclass[11pt,reqno]{amsart}

\newcommand{\dateit}{26 April 2000}

\usepackage{amsmath}

\theoremstyle{plain} 
\newtheorem{theorem}{Theorem}
\newtheorem{proposition}[theorem]{Proposition}
\newtheorem{lemma}[theorem]{Lemma}

\newtheorem{fact}[theorem]{Fact}

\theoremstyle{definition}

\theoremstyle{remark}

\newtheorem{observation}[theorem]{Observation} 

\newtheorem*{Acknowledgement}{Acknowledgement} 

\catcode`@=11 \@mparswitchfalse 

\newcounter{mnotecount}          

\newcommand{\hsp}[1]{\hskip #1 pt}

\newcommand{\e}{\varepsilon}

\newcommand{\w}{\omega}
\newcommand{\X}{\mathfrak X}
\newcommand{\Y}{\mathcal Y}
\newcommand{\N}{\mathbb N}
\newcommand{\R}{\mathbb R}
\newcommand{\lnm}{\left\Vert} 
\newcommand{\rnm}{\right\Vert}

\newcommand{\lsp}{\left[}
\newcommand{\rsp}{\right]}
\newcommand{\lip}{\langle}
\newcommand{\rip}{\rangle_{H_2}}

\newcommand{\auc}{asymptotically uniformly convex}
\newcommand{\cva}{convexifiable}
\newcommand{\auca}{asymptotically uniformly \cva}
\newcommand{\aly}{asymptotically}

\newcommand{\RNP}{Radon-Nikod\'ym property}
\newcommand{\PCP}{point of continuity  property}

\newcommand{\mcx}{\delta_\X}
\newcommand{\amc}{\overline{\delta}}
\newcommand{\amcx}{\overline{\delta}_{\X}}
\newcommand{\fcd}[1]{\mathfrak N \left(#1\right)}
\newcommand{\wfcd}[1]{\mathcal W\left(#1\right)}
\newcommand{\pe}{\widehat{\imath}\:}  

\newcommand{\tr}{\mathcal T}
\newcommand{\ta}{\widetilde{A}}

\newcommand{\diam}{\textrm{diam\,}}
\newcommand{\pperp}{\top}  
\newcommand{\wh}{\widehat}
\newcommand{\tx}{\widetilde{x}\,^*}
\newcommand{\ty}{\widetilde{y}\,^*}

\newcommand{\xy}{x^* +\e y^*}
\newcommand{\tu}{\widetilde{u}\,^*}
\newcommand{\txu}{\tx +  \tu}

\begin{document}

\title[$JT^*$ is asymptotically uniformly convex]{The  
dual of the James Tree space is 
asymptotically uniformly convex}
\author[Maria Girardi]{Maria Girardi}
\address{Dept.~ of Math., University of South Carolina, 
         Columbia, SC   29208, U.S.A.} 
\email{girardi@math.sc.edu}
\date{\dateit}

\begin{abstract} 
The   dual of the James Tree space is  \auc.   
\end{abstract}

\maketitle

\baselineskip 14 pt

\section{Introduction}
\label{s:intro}

In 1950, R.~C.~James \cite{J1}  constructed 
a Banach space which is now called   
the James space. 
This space, along with its many variants 
(such as the James tree space ~ \cite{J2}) 
and their duals and preduals, 
have been a rich source 
for further research and   results 
(both positive ones and counterexamples), 
answering many questions, several of which 
date back to  Banach \cite[1932]{B}.   
See \cite{FG} for a splendid survey of such spaces.

This paper's main result, Theorem~\ref{t:itd},  
shows that  the  dual~$JT^*$ 
of the James tree space ~$JT$  is asymptotically uniformly convex.  
(See Section~\ref{s:dn} for definitions.)

Schachermayer \cite[Theorem~4.1]{S} 
showed that~$JT^*$ 
has the Kadec-Klee property.  
It follows from  
Theorem~\ref{t:itd} of this paper that    
$JT^*$  enjoys   the \textit{uniform} Kadec-Klee property.  
Of course, the same can be said about the 
(unique) predual~$JT_*$ of~$JT$.   
In fact, 
Theorem~\ref{t:pd} shows that the 
modulus of asymptotic convexity of  $JT_*$  
is of power type~3.

Johnson, Lindenstrauss,  Preiss,  and Schechtman~\cite{auc}  
showed that an 
\auc\ space has the \PCP\  
and thus asked 
whether an \auc\ space 
has the \RNP.  
It is well-known that both $JT_*$ and 
$JT^*$ have the  
\PCP\ yet fail  the \RNP. 
It follows from    
Theorem~\ref{t:itd} of this paper that    $JT^*$ 
is an   \auc\  (\textit{dual}) space 
without the \RNP.  
Thus $JT_*$ is a 
\textit{separable} \auc\   space 
without the \RNP. 
To the best of the author's knowledge, these are the 
first  known examples  of   
\auc\ spaces without the \RNP. 
 
\baselineskip 15 pt

\section{Definitions and Notation}
\label{s:dn}

Throughout this paper  $\X$   
denotes  an arbitrary 
(infinite-dimensional real) Banach spaces.   
If $\X$ is a Banach space, then  
$\X^*$  is its dual space,
$B(\X)$ is its  (closed) unit ball,    
$S(\X)$ is its unit sphere, 
$\pe \colon \X \to \X^{**}$ is the natural point-evaluation 
isometric embedding,   $\wh x = \pe (x)$ and 
$\wh{\X} = \pe\left(\X\right)$.    
If $Y$ is a subset of $\X$, then 
$\lsp Y \rsp$  is 
the closed linear span of~$Y$ and
\begin{align*}
\fcd{\X} ~&=~ 
\left\{ \lsp x_i^* \rsp^\pperp_{1\leq i \leq n} 
\colon
x_i^* \in \X^* \textrm{\hsp{4}and\hsp{4}}
n\in\N \right\}\\  
\wfcd{\X^*} ~&=~ 
\left\{ \lsp x_i  \rsp^\perp_{1\leq i \leq n} 
\colon
x_i \in \X \textrm{\hsp{4}and\hsp{4}}
n\in\N \right\} ~. 
\end{align*}     
Thus $\fcd{\X}$ is the collection of 
(norm-closed) finite codimensional subspaces of $\X$ 
while $\wfcd{\X^*}$ is the collection of 
weak-star closed  finite codimensional subspaces of $\X^*$.  
All notation and terminology, not otherwise explained, 
are as in~\cite{DU, LT1, LT2}. 

The \textit{modulus of  convexity} 
$\mcx  \colon [0,2] \to [0,1]$ of~$\X$ is 
$$
\mcx (\e) ~=~ 
\inf \left\{ 1 - \lnm \frac{x+y}{2} \rnm ~\colon~ \hsp{3}
x,y \in S(\X)  \text{\hskip 5 pt and \hskip 5 pt} 
\lnm x - y \rnm \geq \e \right\}            
$$  
and  $\X$ is \textit{uniformly convex (UC)}  if and only if 
$\mcx (\e)>0$ for each $\e\in (0,2]$.  
The \textit{modulus of asymptotic   convexity} 
$\amcx  \colon [0,1] \to [0,1]$ of~$\X$ is 
\begin{equation*} 
\amcx (\e) ~=~ 
\inf_{x\in S(\X)} \hsp{3} 
\sup_{\Y \in\fcd{\X}}      \hsp{3}
\inf_{y\in S(\Y)} \hsp{5} 
\left[ ~\lnm x + \e y \rnm - 1~ \right] 
\end{equation*}
and $\X$ is \textit{\auc\ (AUC)}  if and only if  
$\amcx (\e)>0$ for each $\e$ in ~$(0,1]$.

A space $\X$ has the  \textit{Kadec-Klee (KK) property}  provided 
the relative norm and weak topologies on $B(\X)$    coincide on $S(\X)$. 
A space    
$\X$ has the  \textit{uniform Kadec-Klee (UKK) property}  provided  
 for each $\e >0$ there exists  $\delta >0$ 
such that  
every  $\e$-separated  weakly convergent 
sequence $\{ x_n \}$ in  $B(\X)$  
converges to an element of norm less 
than $1 -\delta$. 

Related to the above  geometric isometric properties 
are the following 
geometric isomorphic  properties.   
\begin{itemize}
\item
$\X$ has the  
\textit{\RNP\ (RNP)} provided  
each  bounded subset of $\X$  has 
non-empty slices of arbitrarily 
small diameter. 
\item 
$\X$ has the 
\textit{point of continuity property\ (PCP)} provided  
each  bounded subset of $\X$  has 
non-empty relatively weakly open subsets of arbitrarily 
small diameter.   
\item
$\X$ has the  
\textit{complete continuity property  (CCP)} provided  
each bounded subset of $\X$ is Bocce dentable.  
\end{itemize}
 
Implications between these various properties 
are summarized in the diagram below.  
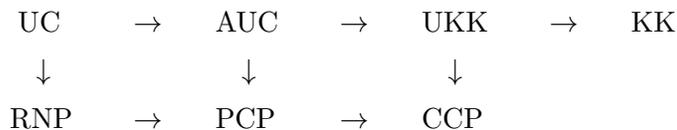
\newcommand{\z}{\hsp{10}}
\begin{alignat*}{7}
&\textrm{\hsp{3}UC\hsp{3}}
\z&\z&\rightarrow
\z&\z&\textrm{AUC}
\z&\z&\rightarrow
\z&\z&\textrm{UKK}
\z&\z&\rightarrow
\z&\z&\textrm{KK}
\\
&\hsp{7}\downarrow
&&&&\hsp{7}\downarrow
&&&&\hsp{7}\downarrow\\
&\textrm{RNP}
\z&\z&\rightarrow
\z&\z&\textrm{PCP}
\z&\z&\rightarrow
\z&\z&\textrm{CCP}
\end{alignat*}

Helpful notation is  
\[
\amcx (\e) ~=~ 
\inf_{x\in S(\X)} \hsp{3} \amcx(\e, x) 
\]
where 
\[
 \amcx(\e, x) ~=~ 
\sup_{\Y \in\fcd{\X}}      \hsp{3}
\inf_{y\in S(\Y)} \hsp{5} 
\left[ ~\lnm x + \e y \rnm - 1~ \right] ~ . 
\] 
Note that, for  each $x\in S(\X)$, 
\[
\amcx (\e, x) ~=~   
\sup_{\Y \in\fcd{\X}}      \hsp{3}
\inf_{\substack{y\in \Y \\ \lnm y \rnm \geq \e }} \hsp{5} 
\left[ ~\lnm x +   y \rnm - 1~ \right]     
 \] 
and so 
$\amcx(\e, x)$ is a non-decreasing function of $\e$. 
Thus  $\amcx$ is non-decreasing Lipschitz functions 
with Lipschitz constant at  most one.  
For any space $\X$ and  $\e\in[0,1]$
\[
\amcx(\e) ~\leq~ \e ~=~ \amc_{\ell_1} (\e) ~;
\]
thus, $\ell_1$ is, in some sense,  the most 
\auc\ space.  

Uniform convexity, the KK property, and the 
UKK property have been extensively studied 
(for example, see \cite{DGZ, LT2}).    
Asymptotic  uniform  convexity has been examined 
explicitly in \cite{auc, M} and implicitly in \cite{GKL, KOS}. 
The RNP,  PCP, and CCP    have also been extensively studied 
(for example, see \cite{DU, GGMS, G1, G2}). 

The  $JT$ space is  construction on  a  (binary) 
\textit{tree} 
\[
\tr ~=~ \mathop{\cup}_{n=0}^\infty \hsp{5} \Delta_n 
\]
where $\Delta_n$ is the $n^{\textrm{th}}$-level of 
the tree; thus, 
\[  
\Delta_0 ~=~ \{ \emptyset \} 
\textrm{\hsp{20}and\hsp{20}}
\Delta_n ~=~ \{ -1, +1   \}^n
\]
for each $n\in\N$.
The \textit{finite tree} $\tr_N$  up through level $N\in \N\cup \{0\}$ is  
\[
\tr_N ~=~ \mathop{\cup}_{n=0}^N \hsp{5} \Delta_n  
\]
The tree $\tr$ is   equipped with its natural (tree)
  ordering: 
if $t_1$ and ~$ t_2  $ are elements of $\tr$,  then  
$
t_1 < t_2
$
provided  one of the follow holds: 
\begin{enumerate}
\item  
 $t_1 = \emptyset$  and $t_2 \not = \emptyset$
\item 
for some  $n,m\in\N$ 
\[ 
t_1 = (\e_1^1, \e_2^1, \ldots, \e_n^1) 
\textrm{\hsp{20}and\hsp{20}} 
t_2 = (\e_1^2, \e_2^2, \ldots, \e_m^2)
\] 
with $n<m$ and $\e_i^1 =\e_i^2$ for each $1\leq  i \leq  n$  ~ . 
\end{enumerate}  
A (finite) \textit{segment} of $\tr$ is a linearly order  subset 
$\{ t_n, t_{n+1}, \ldots, t_{n+k} \}$ of    $\tr$ where  
$t_i\in\Delta_i$ for each $n\leq i\leq n+k$. 
A \textit{branch} of $\tr$ is a linearly order  subset 
$\{ t_0, t_{1},  t_2, \ldots  \}$ of    $\tr$ where  
$t_i\in\Delta_i$ for each $i \in \N \cup \{ 0\}$.

The \textit{James-Tree space} $JT$ is the completion of the   space of 
 finitely supported  functions 
$
 x \colon \tr \to \R
$  
 with respect to the norm
\begin{multline*}
\lnm x \rnm_{JT} ~=~ \\  \sup ~
\left\{ 
\left[ \sum_{i=1}^n \vert \sum_{t\in S_i} x_t \vert^2\right]^{\frac{1}{2}} 
~\colon S_1, S_2, \ldots, S_n \textrm{ are disjoint segments of } \tr 
\right\}  ~. 
\end{multline*} 
By lexicographically ordering $\tr$, 
the sequence $\{ \eta_t \}_{t\in\tr}$ in $JT$, where 
\[
\eta_t(s) ~=~
\begin{cases}
1 &\text{if $t=s$} \\
0 &\text{if $t \not = s$} \ , 
\end{cases}
\]
forms a  monotone 
boundedly complete monotone (Schauder) basis of $JT$ 
with 
biorthogonal functions $\{ \eta_t^* \}_{t\in\tr}$ 
in $JT^*$. 
Thus $\wh{  JT_* } =     \lsp   \eta_t^* \rsp_{t\in\tr} $.

For $N,M \in \N \cup \{ 0 \}$ with $N\leq M$, the restriction maps  
from $JT$ to $JT$ given by 
\begin{align*} 
\pi_N (x) &= \sum_{t\in \Delta_N} \eta_t^* (x) \eta_t \\
\pi_{[N,M]} (x) &= \sum_{t\in \cup_{i=N}^M \Delta_i} \eta_t^* (x) \eta_t \\
\pi_{[N, \w)} (x) 
&= \sum_{t\in \cup_{i=N}^{\infty} \Delta_i} \eta_t^* (x) \eta_t 
\end{align*} 
are each contractive projections 
(by the nature of the norm on $JT$); thus, so are their adjoints. 

Let $\Gamma$ be the set of all branches of $\tr$. 
Then~\cite[Theorem~1]{LS} 
the mapping   
\[ 
\pi_{\infty} \colon JT^* \to \ell_2(\Gamma) 
\]
given by 
\[
\pi_{\infty} (x^*) ~=~ 
\left\{ \lim_{t\in B} x^*(\eta_t) \right\}_{B\in \Gamma} 
\]
is an isometric quotient mapping with kernal $\wh{ JT_*}$. 
Also,  
for each $x^*\in JT^*$, 
\begin{gather*}
\lnm x^* \rnm ~=~ \lim_{N\to\infty} 
\lnm  \pi^*_{[0,N]} \, x^* \rnm 
 \\
\lnm \pi_\infty x^* \rnm ~=~ \lim_{N\to\infty} 
\lnm  \pi^*_{[N, \omega)} \, x^* \rnm~=~
\lim_{N\to\infty} 
\lnm  \pi^*_{N} x^* \rnm  
\end{gather*}
by the weak-star lower semicontinuity of the norm on $JT^*$.

To show that $JT^*$ has the Kadec-Klee property, 
Schachermayer  calculated the below two   quantitative  bounds.  
\begin{fact} 
\label{f:sc1}
{\textrm{\cite[Lemma 3.8]{S}}} 
Let  
\[ 
f_1  \colon  (0,1) \to \left(0, \infty\right) 
\]  
be a continuous strictly increasing function     satisfying  
\(
  f_1(t) < 2^{-10} t^3 
\) 
for each ~$t\in(0,1)$. 
Let $N \in \N$ and $z^* \in JT^*$. If 
\begin{equation*}
\left[   1-f_1(t)   \right] \ \lnm z^* \rnm 
~<~
\lnm \pi^*_{[0,N]} \ z^* \rnm  
\end{equation*} 
then  
\begin{equation*}
\lnm \pi^*_{[N,\w)} \  z^* \rnm 
~<~ 
\lnm \pi^*_{N} \ z^* \rnm  ~+~ t \lnm   z^* \rnm  \ . 
\end{equation*}
\end{fact}
  
\begin{fact} 
\label{f:sc2}
{\textrm{\cite[Lemma 3.11]{S}}} 
Let  
\[ 
f_2  \colon  (0,1) \to \left(0, \infty\right) 
\]  
be a continuous strictly increasing function     satisfying  
\(
  f_2(t) < 2^{-26} t^5 
\) 
for each~ $t\in(0,1)$. 
Let $N \in \N$ and $\e_0\in(0,1)$  and $\tx, \tu  \in JT^*$. If 
\begin{enumerate}
\item[\rm{(2.1)}] 
$\lnm \pi^*_{[N, \omega) } \ \tx \rnm  ~\leq~ 1$ 
\item[\rm{(2.2)}] 
$\lnm \pi^*_{N} \ \tx \rnm ~>~ 1 - f_2(\e_0)$
\item[\rm{(2.3)}] 
$\lnm \pi_{\infty} \ \tx \rnm ~>~ 1 - f_2(\e_0)$
\item[\rm{(2.4)}] 
$\lnm \pi^*_{[N, \omega) } \ (\txu) \rnm  ~\leq~ 1$ 
\item[\rm{(2.5)}] 
$\lnm \pi^*_{N} \, \tu \rnm  ~<~ f_2(\e_0)$
\item[\rm{(2.6)}] 
$\lnm \pi_{\infty}\, \tu \rnm  ~<~ f_2(\e_0)$ . 
\end{enumerate} 
Then 
\begin{enumerate}
\item[\rm{(2.7)}]
$\lnm \pi^*_{[N, \omega) } \,  \tu \rnm  ~<~  \e_0 $ .
\end{enumerate}
\end{fact}

\section{Results}

Theorem \ref{t:pd} shows that the 
modulus of asymptotic convexity of  $JT_*$  
is of power type 3.  Its proof uses  Fact~\ref{f:sc1}.

\begin{theorem}
\label{t:pd}
There exists a positive constant $k$ so that 
\[
 \amc_{JT_*}(\e) ~\geq~ k \e^3 
\]
for each $\e\in(0,1]$.  Thus $JT_*$ is \auc. 
\end{theorem} 

\begin{proof} 
Fix $c\in\left(0, 2^{-10}\right)$ and find $k$ so that 
\begin{equation}
\label{eq:pdk}
0 ~<~ k(1+k)^2 ~\leq~ c ~ .  
\end{equation}
Fix $\e\in(0,1)$ and a finitely supported $x_* \in S\left(JT_*\right)$. 
It suffices to show that 
\begin{equation}
\label{eq:pdwts}
 \amc_{JT_*}(\e, x_*) ~\geq~ k \e^3  ~ . 
\end{equation}

Find $N\in \N$ so that 
\[
  \pi^*_{[0, N-1]}\, \wh x_*  ~=~ \wh x_* 
\]
and let 
\[
\Y ~=~ \lsp \eta_t \rsp^{\pperp}_{t\in \tr_{N}} \ . 
\]
Fix $y_* \in S(\Y)$.  

Assume that 
\begin{equation*}
 \lnm   x_* + \e y_* \rnm ~-~ 1 ~<~ k \, \e^3 \ . 
\end{equation*}
Then 
\[
\left[ 1 -  \frac{k\e^3}{1+k\e^3} \right] \,
 \lnm \wh x_* + \e \wh y_* \rnm
~<~  1  
~=~ \lnm \pi^*_{[0,N]} \, \left(  \wh x_* + \e \wh y_*   \right) \rnm ~ . 
\]
Thus by Fact~\ref{f:sc1}, with $f_1(t) = ct^3$, 
\[
\lnm \pi^*_{[N,\w)} \  \left( \wh x_* + \e \wh y_*  \right) \rnm 
~<~ 
\lnm \pi^*_{N} \ \left( \wh x_* + \e \wh y_*  \right) \rnm  
~+~ f_1^{-1} \left(\frac{k\e^3}{1+k\e^3} \right) \lnm  
\,  \left( \wh x_* + \e \wh y_*  \right) \rnm   
\]
and so 
\begin{equation}
\label{eq:pdc}
\e  ~<~
\left[ 1 +   k\e^3  \right] \,  f_1^{-1} \left(\frac{k\e^3}{1+k\e^3} \right)
  ~. 
\end{equation}
But inequality \eqref{eq:pdc} is equivalent to 
\[
c^{1/3} ~<~ k^{1/3} \left( 1+k\e^3\right)^{2/3} ~ , 
\] 
which contradicts \eqref{eq:pdk}.  
Thus $\lnm   x_* + \e y_* \rnm ~-~ 1 ~\geq~ k \, \e^3 $ 
and so ~\eqref{eq:pdwts} holds.  
\end{proof}

A modification of the proof of Theorem~\ref{t:pd} 
shows   that, 
for  each  $\e\in(0,1)$, 
the \(
 \amc_{JT^*}(\e, x^*)  
\) 
stays uniformly bounded below from zero
for $x^*\in S(JT^*)$ whose $\lnm \pi_\infty \, x^* \rnm$ 
is small.  
Recall that if $x_* \in JT_*$ 
then  $\lnm \pi_\infty  \, \wh x_* \rnm = 0$.

\begin{lemma} 
\label{l:it} 
For each $\e\in (0,1)$  there exists $\eta = \eta(\e)>0$ so that    
\[ 
\inf_{\substack{ x^*\in S(JT^*) \\
        \lnm \pi_\infty x^*\rnm \leq \eta} } \hsp{6} 
\sup_{\Y \in\wfcd{JT^*}}      \hsp{5}
\inf_{y^*\in S(\Y)} \hsp{5} 
\left[ ~\lnm  x^* + \e y^* \rnm - 1~ \right] 
~>~0 ~. 
\]
\end{lemma}

\begin{proof}
Fix $\e \in (0,1)$. Keeping with the notation in 
Fact~\ref{f:sc1}, find $\delta , \eta_2 >0$ so that 
\begin{gather*}
4\eta_2 ~+~ \frac{\delta}{1-f_1(\delta)} 
~<~ \e   ~. 
\end{gather*}   
Fix $x^*\in S(JT^*)$ 
with 
\[
\lnm \pi_\infty x^* \rnm ~\equiv~b \leq \eta_2 ~. 
\]  
It suffices to show that 
\begin{equation} 
\label{eq:wts}  
\sup_{\Y \in\wfcd{JT^*}}      \hsp{3}
\inf_{y^*\in S(\Y)} \hsp{5} 
\lnm   x^* + \e y^* \rnm ~\geq~ \frac{1}{1-f_1(\delta)} ~~.
\end{equation}

Fix $\eta_1 \in \left(0,1\right)$.  
Find $N \in \N$ so that 
\[
1-\eta_1 \leq \lnm \pi^*_{[0,N]} \,   x^* \rnm 
\textrm{\hsp{20}and\hsp{20}}
\lnm \pi^*_{[N,\w)}\,   x^* \rnm  ~<~ b+\eta_2  
\]
and let 
\[
\Y ~=~ \lsp \eta_t \rsp^{\perp}_{t\in \tr_{N}} \ . 
\]
Fix $y^* \in S(\Y)$.  

Assume that 
\begin{equation*}
 \lnm   x^* + \e y^* \rnm ~<~ \frac{1-\eta_1}{1-f_1(\delta)} \ . 
\end{equation*}
Then 
\[
\left[ 1 - f_1(\delta) \right] \, \lnm  x^* + \e y^* \rnm
~<~ 
\lnm \pi^*_{[0,N]} \,  x^* \rnm 
~=~ \lnm \pi^*_{[0,N]} \, \left(  x^* + \e y^* \right) \rnm ~ . 
\]
Thus by Fact~\ref{f:sc1}
\[
\lnm \pi^*_{[N,\w)} \  \left(  x^* + \e y^* \right) \rnm 
~<~ 
\lnm \pi^*_{N} \ \left(  x^* + \e y^* \right) \rnm  
~+~ \delta \lnm   \left(  x^* + \e y^* \right) \rnm   
\]
and so 
\[
\e - ( b+\eta_2 )
~<~
( b+\eta_2 ) +  \frac{\delta}{1-f_1(\delta)} ~. 
\] 
But $b\leq \eta_2$ and so 
\[
\e  
~<~
4\eta_2   +  \frac{\delta}{1-f_1(\delta)} ~. 
\]  
A contradiction, thus 
\begin{equation*}
 \lnm  x^* + \e y^* \rnm ~\geq~ \frac{1-\eta_1}{1-f_1(\delta)} \ . 
\end{equation*} 
Since $\eta_1>0$ was arbitrary,  inequality~\eqref{eq:wts} holds.   
\end{proof}

Thus  to show that  $JT^*$ is 
\auc, one just needs to examine  
 \(
 \amc_{JT^*}(\e, x^*)  
\) 
for $x^*\in S(JT^*)$ whose $\lnm \pi_\infty x^* \rnm$ 
is  not small.  
Fact~\ref{f:sc2} is used for this case. 
 
\begin{theorem} 
\label{t:itd} 
  $JT^*$ is \auc. 
\end{theorem}

\begin{proof}
Fix $\e\in(0,1)$ and let   
 $\e_0 =  {\e}/{4}$. Let $f_1\colon(0,1)\to(0,2^{-12})$ be given 
by $f_1(t)= 2^{-12} t^3$ and $f_2$ be a function satisfying the 
hypothesis in Fact~\ref{f:sc2}.   
Find $\delta , \eta_2>0$ so that 
\begin{gather*}
4\eta_2 + \frac{\delta}{1-f_1(\delta)} 
~<~ \e   ~. 
\end{gather*}  
Next find $\gamma_i>0$ and $\tau>1$ so that 
\begin{gather}
\gamma_3 ~<~ \gamma_2 ~<~ \frac{1}{2} 
\label{eq:g}
\\
\tau \leq \frac{(1-\gamma_1)\,(1-\gamma_2)}{1-f_2(\e_0)} 
\label{eq:d1}
\\
\tau < \frac{ 1-\gamma_2 }{\sqrt{1-f^2_2(\e_0)}} 
\label{eq:d2}
\\
\tau \leq \frac{\eta_2^3 \gamma_3^3}{2^{15}(1-\gamma_2)^3}
   ~-~ \gamma_4 ~+~ 1 
\label{eq:d3}
\\
\frac{\tau - 1 + \gamma_4}{\tau} ~<~ f_1(1) 
\\
\tau ~\leq ~\frac{1}{1-f_1(\delta)} 
\label{eq:pd} ~.
\end{gather}
 
Fix $x^*\in S(JT^*)$. 
It suffices to show that 
\begin{equation}
\label{eq:wtsd}
\sup_{\Y \in\fcd{JT^*}}      \hsp{3}
\inf_{y\in S(\Y)} \hsp{5} 
\lnm x^* + \e y^* \rnm  ~\geq~ \tau ~. 
\end{equation}
Let 
\[
\lnm \pi_\infty x^* \rnm ~\equiv~b   ~. 
\]   
If $b \leq \eta_2$,   
then by the proof of 
Lemma~\ref{l:it} and \eqref{eq:pd}, inequality~\eqref{eq:wtsd} holds. 
So let $b > \eta_2$.  
Find $N\in \N$ so that 
\begin{gather}
(1-\gamma_1)\,  b 
~<~ 
\lnm \pi^*_N \, x^* \rnm
~\leq~
\lnm \pi^*_{[N, \omega)} \, x^* \rnm 
~<~ 
b \, \left( \frac{1-\gamma_3}{1-\gamma_2} \right)
~<~ 
\frac{b}{1-\gamma_2}
\label{eq:tail}
\\
1-\gamma_4~<~ \lnm \pi^*_{[0,N]} \, x^* \rnm ~. 
\label{eq:start}
\end{gather} 
Let $g_{x^*}\in JT^{**}$ be the functional given by 
\[
g_{x^*}(z^*) = \lip \pi_\infty z^* ,  \pi_\infty x^* \rip 
\]
where the inner product in the natural inner product 
on~$\ell_2(\Gamma)$.   
Let 
\[
\Y ~=~ \lsp \eta_t \rsp^{\perp}_{t\in \tr_{N}} 
~\cap \lsp g_{x^*} \rsp^{\pperp}   
\]
and fix $y^* \in S(\Y)$. 

Assume that 
\begin{equation}
\label{eq:c}
\lnm x^* + \e y^* \rnm ~<~ \tau \ . 
\end{equation} 
It suffices to find a contradiction to~\eqref{eq:c}.   
Towards this, let 
\[
\tx ~=~ \frac{1-\gamma_2}{\tau b} ~~  x^* 
\textrm{\hsp{20}and\hsp{20}}  
\ty ~=~ \frac{1-\gamma_2}{\tau b} ~~  y^* ~.
\] 
It suffices  to show (keeping with 
the same notation but with $\tu = \e \ty$)   that
conditions (2.1) through (2.6) of Fact~\ref{f:sc2} hold; 
for then condition (2.7) holds and so by~\eqref{eq:g}
\[
\e_0 ~>~ 
\lnm \pi^*_{[N, \omega) } \, \e\ty \rnm 
~=~ 
\frac{1-\gamma_2}{\tau b} \ \e  
~\geq~
\frac{\e}{4}  
~=~ \e_0~. 
\]

Condition~(2.1) follows from  \eqref{eq:tail}  since  
\[
\lnm \pi^*_{[N, \omega) } \ \tx \rnm  
~\leq~
\frac{1-\gamma_2}{\tau b} \, \frac{b}{1-\gamma_2} 
~\leq~1 \ . 
\] 
Condition (2.2) follows from \eqref{eq:tail} and \eqref{eq:d1} since 
\[
\lnm \pi^*_{N} \ \tx \rnm 
~>~
\frac{1-\gamma_2}{\tau b} \, (1-\gamma_1)\,  b 
~=~ 
\frac{(1-\gamma_1) \, (1-\gamma_2)}{\tau} 
~\geq~ 1-f_2(\e_0) ~ .
\]
Towards 
condition (2.3), note that by    \eqref{eq:d2}  
\begin{equation}
\label{eq:l2}
\lnm \pi_{\infty} \ \tx \rnm 
~=~ 
\frac{1-\gamma_2}{\tau b} \  b 
~=~
\frac{1-\gamma_2}{\tau}
~>~
\sqrt{1-f^2_2(\e_0)}
\end{equation}
and so 
\[
\lnm \pi_{\infty} \ \tx \rnm 
~>~
1 - f_2(\e_0) ~ . 
\]  

Towards condition (2.4), note that by \eqref{eq:c} 
and ~\eqref{eq:start} 
\[
\lnm \xy \rnm ~<~ 
\frac{\tau}{1-\gamma_4} \ \lnm \pi^*_{[0,N]}\, \left( \xy \right) \rnm ~.
\]
Thus by Fact~\ref{f:sc1} and (9) 
\begin{align*}  
&\lnm  \pi^*_{[N, \omega)}        \left(\xy\right) \rnm \\
~&\hsp{40}<~  
\lnm \pi^*_N \left(\xy\right) \rnm
~+~ f_1^{-1} \left(\frac{\tau-1+\gamma_4}{\tau} \right) \ 
 \lnm   \left(\xy\right) \rnm  \\
~&\hsp{40}\leq~  b \ \frac{1-\gamma_3}{1-\gamma_2} ~+~ \tau 2^4 
\left( \frac{\tau-1+\gamma_4}{\tau} \right)^{1/3} \ . 
\end{align*}
Thus condition (2.4) holds provided 
\[
 b \ \frac{1-\gamma_3}{1-\gamma_2} ~+~ \tau 2^4 
\left( \frac{\tau-1+\gamma_4}{\tau} \right)^{1/3}
~\leq~ 
\frac{\tau b}{1- \gamma_2} \ ,
\]
or equivalently 
\[ 
\tau^{2/3}\ \left( \tau-1+\gamma_4\right)^{1/3} 
~\leq~
\frac{b\left(\tau - 1 +\gamma_3\right)}{2^4\left(1-\gamma_2\right)} \ . 
\]
But by \eqref{eq:d3} and that $b>\eta_2$ 
\begin{alignat*}{2}
\tau^{2/3}\ \left( \tau-1+\gamma_4\right)^{1/3} 
~&~\leq~
2 \ \left( \tau-1+\gamma_4\right)^{1/3}
&&~\leq~
\frac{2\eta_2 \gamma_3}{2^{5}(1-\gamma_2)}\\
&~\leq~
\frac{  b \gamma_3}{2^{4}(1-\gamma_2)}
&&~\leq~
\frac{b\left(\tau - 1 +\gamma_3\right)}{2^4\left(1-\gamma_2\right)} \ . 
\end{alignat*}
Thus condition (2.4) holds.

Condition  (2.5)  follows from the fact   
that  $y^* \in \lsp \eta_t \rsp^{\perp}_{t\in \tr_{N}}$. 
Towards  condition~(2.6),    
since $y^* \in \lsp g_{x^*} \rsp^{\pperp}$,  
the vectors  $\pi_\infty   \ty$ and $\pi_\infty \tx$   
are orthogonal   in ~$\ell_2(\Gamma)$ and so  
\[
\lnm \pi_{\infty}  \e \ty \rnm^2 
~=~ 
\lnm \pi_{\infty} (\tx + \e \ty) \rnm^2 
~-~ 
\lnm \pi_{\infty}  \tx \rnm^2 \ ; 
\]
but $\pi_\infty = \pi_\infty \pi^*_{[N, \omega)}$ 
and so by~ condition~(2.4) and ~\eqref{eq:l2}   
\begin{align*}
\lnm \pi_{\infty}  \e \ty \rnm^2 
~&\leq~
\lnm \pi^*_{[N, \omega) } (\tx + \e \ty) \rnm^2  
~-~ 
\lnm \pi_{\infty}  \tx \rnm^2 \\
~&<~ 
1 ~-~ \left[ 1 - f^2_2(\e_0) \right] 
~=~ f^2_2(\e_0)   . 
\end{align*}
Thus condition (2.6). 
 \end{proof}

The proof in \cite{auc} 
that an \auc\  space has the 
PCP show that if $\amcx(\e) >0$ for each $\e\in (0,1]$ 
then $\X$ has the PCP.  
A bit more can be said.

\begin{proposition}
\label{p:pcp} 
If  $\amcx\left(\frac{1}{2}\right) >0$  then $\X$ has the PCP.
\end{proposition}

\noindent
The proof of Proposition~\ref{p:pcp} uses 
the following (essentially known) lemma. 

\begin{lemma} 
\label{l:pcp}  
Let  $\X$ be a space without the PCP and $0<\e<1$.
Then there is a closed subset~ $A$ of $\X$ so that 
\begin{enumerate} 
\item[\rm{(1)}] 
each (nonempty) relatively weakly open subset of $A$ 
has diameter larger than $ 1 - \e$ 
\item[\rm{(2)}] 
$\sup \{ \, \lnm a \rnm \colon a\in A \} ~=~ 1$. 
\end{enumerate}
\end{lemma} 

\begin{proof}[Proof of Lemma \ref{l:pcp}]
Let $\X$ fail  the PCP and $0<\e<1$.  
By a standard argument (e.g., see   \cite[Prop.~4.10]{SSW}), 
there is a closed subset $\ta$ of $\X$  
of diameter one such that 
each (nonempty) relatively weakly open subset of $\ta$ 
has diameter larger than $1-\e$.  
Without loss of generality  $0\in \ta$ 
(just consider a translate of~ $\ta$). 
Let 
\[
b = \sup\, \{\, \lnm x \rnm \colon x\in \ta \} 
\textrm{\hsp{20}and\hsp{20}} A = \frac{\ta}{b} \, ~.
\]
Note that $0 < b \leq 1$.  
If $V$ is (nonempty) relatively weakly open subset of~ $A$, 
then $bV$ is a   relatively weakly open subset of $\ta$ 
and so 
\[
\diam V = \frac{1}{b}\  \diam bV > 1- \e ~. 
\]
Thus $A$ does the job. 
\end{proof}

\begin{proof}[Proof of Proposition \ref{p:pcp}] 
Let $\X$ be a Banach space without the PCP. 
Fix ~$ t\in \left(0, \frac{1}{2}\right)$ and $\delta\in(0,t)$.  
It suffices to show that $\amcx(t) \leq 2 \delta$.

Find a subset $A$ of $\X$ which satisfies the conditions 
of Lemma~\ref{l:pcp} with $\e=1-2t$ and find 
$a\in A$ so that 
\[ 
\lnm \frac{a}{\lnm a\rnm} ~-~ a \rnm ~<~\delta ~. 
\]
Let $\Y \in \fcd{\X}$.  It suffices to show that 
\[
\inf_{\substack{y\in \Y \\ \lnm y \rnm \geq t }}   \hsp{5} 
\left[~ \lnm  \frac{a}{\lnm a\rnm} +   y \rnm - 1 ~\right]
~\leq~ 2 \, \delta ~.
\]

By condition (1) of Lemma~\ref{l:pcp}  there exists $x\in A$ so that 
$\lnm x-a \rnm \geq t $ and $x-a$ is 
\textit{almost} in $\Y$; thus,   
by a standard  perturbation argument 
(e.g., see~\cite[Lemma~2]{GJ}) there exists 
${y} \in \Y$ so that 
\[
\lnm  {y} \rnm \geq t 
\textrm{\hsp{10}and\hsp{10}} 
\lnm   {y} - \left( x-a\right)  \rnm ~<~ \delta ~. 
\]
Thus 
\[
\lnm  \frac{a}{\lnm a\rnm} +   y \rnm 
~\leq~
\lnm  \frac{a}{\lnm a\rnm} - a \rnm 
~+~
\lnm  {y} - x + a    \rnm 
~+~
\lnm    x  \rnm  
<~ 1 + 2\delta ~.
\]
Thus $\amcx (\frac{1}{2}) = 0$.
\end{proof}

The observation below formalizes  
an essentially known fact, which to the best of the 
author's knowledge, has not appeared in print as such. 
Recall that the \textit{modulus of asymptotic   smoothness}  
$\overline{\rho}_{\X}  \colon [0,1] \to [0,1]$ of~$\X$ is 
\[ 
\overline{\rho}_{\X} (\e) ~=~ 
\sup_{x\in S(\X)} \hsp{3} 
\inf_{\Y \in\fcd{\X}}      \hsp{3}
\sup_{y\in S(\Y)} \hsp{5} 
\left[ ~\lnm x + \e y \rnm - 1~ \right] 
\]
and $\X$ is \textit{asymptotically  uniformly   smooth}  if and only if 
$\lim_{\e \to 0^+}  {\rho_{\X} (\e) }/{\e} = 0$.  
Also,  $L_p(\X)$ is the 
Lebesgue-Bochner space 
of strongly measurable $\X$-valued functions 
defined on a separable non-atomic probability space, 
equipped with is usual norm.

\begin{observation}
Let $1<p<\infty$.  For a Banach space $\X$, 
the following are equivalent. 
\begin{enumerate} 
\item[\rm{(1)}] 
$\X$ is  uniformly \cva.
\item[\rm{(2)}] 
$L_p(\X)$ is  uniformly \cva.
\item[\rm{(3)}] 
$L_p(\X)$ is \auca. 
\item[\rm{(4)}] 
$L_p(\X)$ admits an equivalent UKK norm. 
\item[\rm{(5)}] 
$L_p(\X)$ is \aly\ uniformly smoothable.  
\end{enumerate}  
\end{observation}

\begin{proof} 
Let $1<p<\infty$ and $\X$ be a Banach space.   

That (1) though (4)  are equivalent and that~ (2) implies ~(5) 
follows easily from 
the below known facts about a Banach space $\Y$. 
\begin{enumerate}
\item[\rm{(i)}]  
$\Y$ is uniformly convex if and only if  
$L_p(\Y)$ is 
\cite{Mc}.   
\item[\rm{(ii)}] 
$\Y$ is uniformly \cva \ 
if and only if $L_p(\Y)$  admits an equivalent UKK norm 
\cite[Theorem 4]{DGK}.
\item[\rm{(iii)}]  
$\Y$ is uniformly \cva\ if and only if 
$\Y$ is uniformly smoothable 
(cf.~\cite[page 144]{DU}).
\end{enumerate}

Towards showing that (5) implies (1), 
let   $L_p(\X)$ be \aly\ uniformly smoothable 
and $\X_0$ be a separable subspace of $\X$. 
It suffices to show that $\X_0$ is uniformly \cva\ 
(cf.~\cite[Remark~IV.4.4]{DGZ}).

It follows from   \cite[Proposition~2.6]{GKL} 
 that if $\Y$ is separable, then 
$\Y$ is \aly\ uniformly smooth if and only if 
$\Y^*$ has the UKK$^*$ property.  
Thus $\left[L_p(\X_0)\right]^*$  admits an equivalent 
UKK$^*$ norm. 
But   $\ell_1$ cannot embed into  
$L_p(\X_0)$ since ~$L_p(\X_0)$ is    \aly\ uniformly smoothable    
and so $\left[L_p(\X_0)\right]^*$  is 
\aly\ weak$^*$ uniformly \cva\ and so is also 
\aly\  uniformly \cva.  
Thus $L_q(\X_0^*)$ is \auca\ 
where $1/p + 1/q = 1$.  
From (3) implies (1) it follows that 
$\X_0^*$ is uniformly \cva\ and so so is $\X_0$.  
\end{proof}

\begin{Acknowledgement} 
The author thanks  William B.~Johnson 
and Thomas Schlumprecht 
for their fruitful discussions  on 
asymptotic  uniform  convexity  
at the NSF Workshop in Linear Analysis and Probability, 
Texas A\&M University, during the Summer of 1999.  
\end{Acknowledgement}

\end{document}